# Development and operational deployment of a staff rostering system using mathematical optimisation in a municipal water utility

Takaaki Kobayashi[1], Satoshi Matsumura[2] and Satoko Moriguchi[3],*

[1] TransRecog Co., Ltd., Tokyo, Japan

[2] Water Purification Group, Water Purification Office, Suita City Waterworks Bureau, Suita, Osaka, Japan

[3] Faculty of Economics and Business Administration and Graduate School of Management, Tokyo Metropolitan University, Hachioji, Tokyo, Japan

* Corresponding author: Satoko Moriguchi (satoko5@tmu.ac.jp); ORCID: 0000-0003-0105-7167

Takaaki Kobayashi: takobaya@transrecog.com; ORCID: 0009-0006-8492-5677

Satoshi Matsumura: w-jyosui@city.suita.osaka.jp; ORCID: 0009-0000-6083-6588

## ABSTRACT

We developed a staff rostering system for complex employment and staffing requirements at a municipal water utility and deployed it in routine operations at the Water Purification Office of Suita City Waterworks Bureau, Japan, in May 2026. We estimated that preparing a monthly roster took approximately 2–3 h after deployment, compared with 12–18 h previously. The system represents roster states, effective staffing and office assignments using binary and auxiliary integer variables and optimises the entire planning horizon. When the normal model yields no feasible solution, a second stage attempts to generate a reference roster by penalising weekly day-off shortages while retaining the other hard constraints. Outputs identify shortages and meeting-date changes applied before solving. Computational experiments used operational inputs for 14 staff over 245 days. The baseline was solved to optimality in approximately 177 s; a modified meeting-day scenario produced a reference roster with five staff-days of shortage. Independent monthly solving without history produced 16 staff–week cases of excess days off across two month-boundary weeks. This case demonstrates the practical value of optimisation-based decision support for staff rostering in municipal water utilities by showing how mathematical optimisation can be integrated into roster preparation and checking workflows with reduced preparation time.

Keywords: CP-SAT; decision support; integer optimisation; staff rostering; water utilities

## HIGHLIGHTS

1. Optimisation supports routine staff rostering at a municipal water utility.
2. Estimated monthly roster preparation time fell from 12–18 h to about 2–3 h.
3. A single model evaluates weekly rules across month boundaries.
4. Conditional second-stage solving identifies weekly day-off shortages.
5. Checking sheets show shortages and applied meeting-date changes.

## GRAPHICAL ABSTRACT

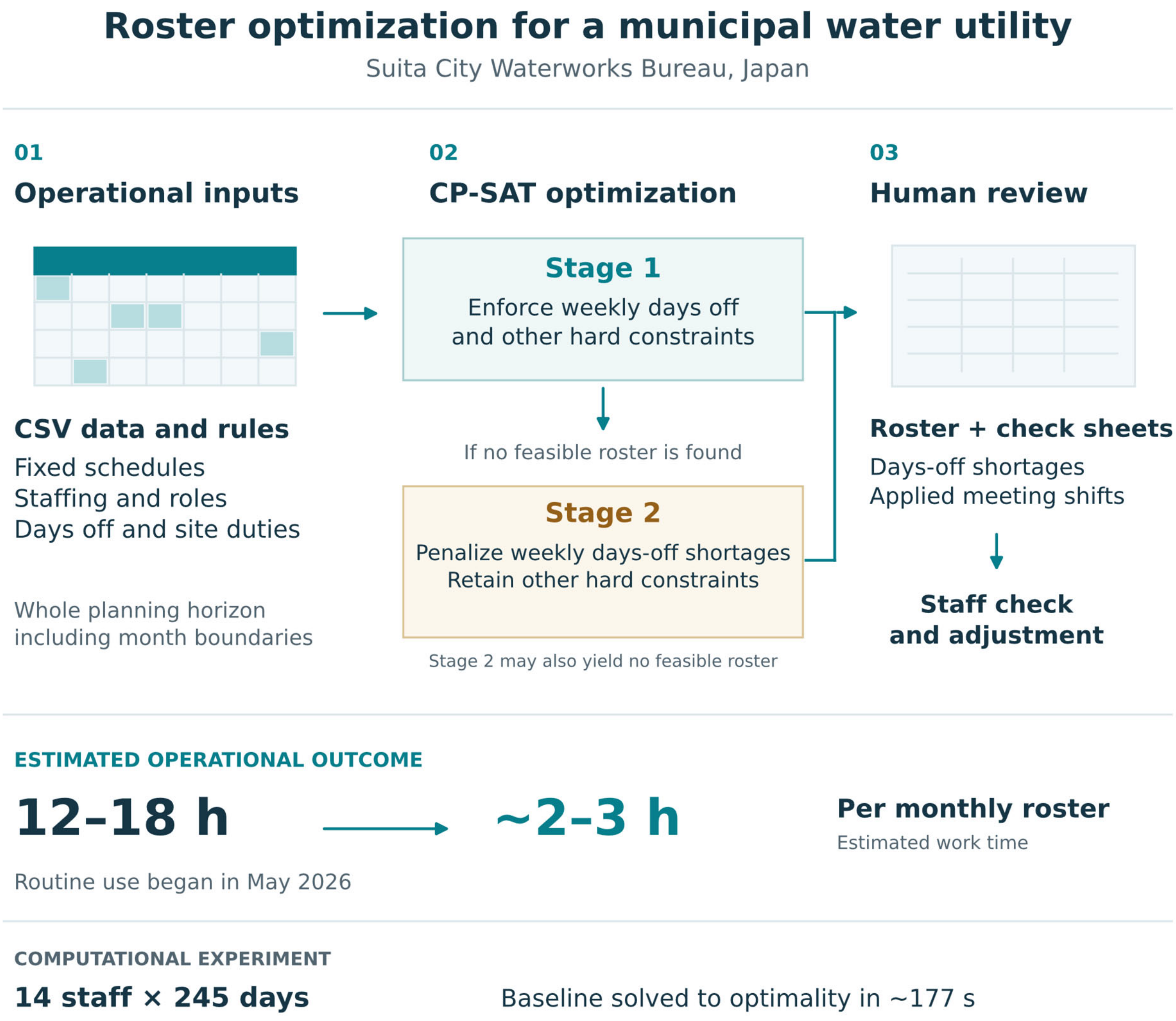

# 1 INTRODUCTION

## 1.1 Background

Water utilities must maintain staffing for daily operations and responses to leaks and equipment failures. Rostering in the department studied requires simultaneous consideration of fixed schedules, days off by employment category, consecutive working days, office assignments and managerial coverage. These requirements interact, and weekly rules cross month boundaries, making manual preparation an iterative process of adjustment and checking.

This study reports how operational work rules were formulated as a mathematical optimisation model and how roster generation and checking information were incorporated into routine municipal operations. We describe deployment and operational effects at the Water Purification Office of Suita City Waterworks Bureau and evaluate the supporting method through computational experiments. The following review places the work within decision support for water utilities, personnel rostering, and workforce and maintenance planning for public infrastructure.

## 1.2 Operational management and decision support in water utilities

Decision support in water utilities encompasses operations, asset management, performance assessment and digital transformation. Miltell et al. (2025) developed a digital transformation roadmap addressing organisational, cultural and interoperability issues. Bettin (2023) examined a business intelligence platform integrating operational and business data; Dadiala & Sinha (2025) considered data governance and modelling; and Lumley et al. (2024) addressed integrated urban water management using digital twins. These studies emphasise the connection between computational tools, operational data and decision processes.

Abe & Mizutani (2016) reported internal auditing, certification and maturity assessment based on ISO 55001 at Sendai's wastewater utility, while Abdelghany & Abdel-Monem (2019) used a balanced scorecard for utility performance assessment. Both make operational conditions and evaluation criteria explicit to support improvement.

Studies of utility performance include Berg's (2013) review of related books, Berg & Marques' (2011) survey of quantitative research, the international comparison by da Cruz et al. (2012), and the analysis of governance, strategy and efficiency by Romano et al. (2018). Beuken et al. (2020) reviewed asset management research, and Morley et al. (2016) compared intervention strategies through a decision support system based on urban water metabolism modelling.

Whereas these studies principally concern assets, networks, operations or utility performance, this study focuses on daily staff allocation. It treats compliance with roster requirements and checking by the responsible staff member as parts of the same operational workflow.

## 1.3 Mathematical optimisation for staff scheduling and rostering

Staff scheduling and rostering have an extensive research literature. Ernst et al. (2004) reviewed industries, models and solution methods, while Cheang et al. (2003), Burke et al. (2004) and Van den Bergh et al. (2013) surveyed problems involving varied work rules and institutional settings. Böðvarsdóttir & Stidsen (2026) found that weighted single-objective formulations and goal programming are widely used in practical multi-objective personnel rostering.

Relevant modelling work includes shifts for employees of different seniority levels (Topaloglu, 2009), representations of assignment history from previous periods (Kingston, 2021), and graph-based formulations for complex work rules (Lai et al., 2020). These provide context for managerial coverage and the treatment of planning-period boundaries.

Solution approaches include heuristics and combinations of integer programming and neighbourhood search (Burke et al., 2008, 2010), a two-phase method (Valouxis et al., 2012), integer programming (Santos et al., 2016), hybrid integer and constraint programming (Rahimian et al., 2017), and a mixed integer programming system for operational use (Böðvarsdóttir et al., 2022). Valouxis et al. (2012) separate workload allocation from individual assignments. Our second stage instead relaxes weekly day-off requirements when the normal model yields no feasible solution, giving the stages a different purpose and structure.

## 1.4 Public infrastructure and field service scheduling

In related fields, Castillo-Salazar et al. (2016) surveyed workforce scheduling and routing, and Borenstein et al. (2010) examined the effects of partitioning dynamic workforce scheduling problems. Xu & Chiu (2001), Kovacs et al. (2012) and Zamorano & Stolletz (2017) studied technician assignment and routing with skills, time windows or multiple planning days. The branch-office assignments in our system likewise involve coordinating staff across locations.

For water infrastructure, Salman et al. (2013) considered distribution network rehabilitation, Osman et al. (2017) studied multi-objective repair-crew scheduling, Draude et al. (2022) examined sewer maintenance scheduling with crew travel, and Zaman et al. (2017) addressed productivity in drainage maintenance. Our study is related to these workforce and task planning problems but concerns the daily assignment of work and days off.

### 1.5 Positioning and contribution

The principal contribution is the implementation of a municipal water department's rostering requirements in a mathematical optimisation model and the deployment of the resulting system in routine operations. Operational use began at the Water Purification Office of Suita City Waterworks Bureau in May 2026. Operational estimates indicate a reduction in preparation time per monthly roster from 12–18 h to approximately 2–3 h.

The study connects water utility decision support with personnel rostering by showing how work rules can be integrated into input preparation, solving, checking and output. The computational experiments evaluate the formulation and exception handling that support this operational application.

Three technical features support this workflow. The model distinguishes roster states from effective staffing and jointly represents employment-specific days off, office and managerial assignments, fixed schedules and period boundaries. Conditional two-stage solving provides normal roster generation followed, when needed, by an attempt to generate a reference roster with shortages. Outputs identify shortages and meeting-date changes already applied, providing information for staff review.

## 2 STUDY SITE AND PROBLEM DESCRIPTION

### 2.1 Water utility and service continuity

The study concerns staff rostering at the Water Purification Office of Suita City Waterworks Bureau in Osaka Prefecture, Japan. Computational experiments used operational inputs for the 14 staff included in the roster studied. Table 1 summarises the utility as a whole, based on its annual report (Suita City Waterworks Bureau, 2025).

The roster assigns daily work, days off and office locations. It is an operational plan for maintaining the headcount and managerial coverage needed for routine operations and responses to leaks and equipment failures.

The scope is the directly staffed daytime service from 09:00 to 18:00, including Saturdays and Sundays. Night-time and early-morning water-supply monitoring and emergency response are provided by an external contractor and are outside this study's scope.

***Table 1. Overview of the water utility (Suita City Waterworks Bureau).***

| Item | Description |
|---|---|
| Water utility | Suita City Waterworks Bureau, Osaka, Japan |
| Population served | 383,913 (FY2024) |
| Annual distributed water volume | 41,853,560 $m^3$ (FY2024) |
| Pipeline length | 746,530 m (as of 31 March 2025) |
| Purification facilities | Izumi and Katayama purification plants |

Source: Suita City Waterworks Bureau, Waterworks Business Annual Report, FY2025 edition (Suita City Waterworks Bureau, 2025), pp. 33, 44, 50 and 52. These are utility-wide figures. Population and distributed water volume refer to FY2024; pipeline length is as of 31 March 2025.

## 2.2 Rostering requirements

Table 2 summarises the main requirements, which fall into five groups.

First, ordinary weekdays require at least the specified effective headcount, whereas weekends, public holidays and office-closed days require exactly that headcount. Business trips, training and leave do not count towards effective staffing. Staff are trained and deployed to perform the department's range of tasks; the model therefore does not assign tasks by skill but separately represents managerial grades needed for judgement and supervision.

Second, work takes place at the main office and a branch office. A specified number of staff must work at the branch each day. Because an individual's location may vary by day, work, days off and office assignments are determined jointly.

Third, at least one section chief or higher must be present each day, and at least one must remain at the main office rather than all being assigned to the branch.

Fourth, the model reflects weekly work and day-off requirements for regular, contract and part-time staff, together with fixed days off, holidays, leave and training. Group-meeting days normally count as effective work, but fixed weekly days off and fixed holiday days off take precedence. A meeting falling on a non-business day is moved to the next business day. Consecutive working days and consecutive ordinary days off are also limited. Section 3 defines rule precedence and exceptions.

Fifth, weekly day-off and consecutive-day rules are evaluated across month boundaries. If a week spans two months, a monthly solution approach must carry assignment information across the boundary to count the whole week.

Mandatory conditions also include balancing monthly branch-office assignments and limiting full-horizon differences in holiday and weekend days off among regular and contract staff. Monthly imbalance, branch duty repeated after short intervals, and a target of two managers are evaluated through a weighted objective. Weekly days off are an equality within the normal model's defined scope; exceptions and second-stage shortages are treated separately.

***Table 2. Conditions to be considered in roster creation.***

| Category | Main conditions | Treatment |
|---|---|---|
| Effective staffing | At least 6 staff on ordinary weekdays; exactly 6 on weekends, public holidays and office-closed days. | Hard |
| Branch-office duty | Exactly 2 staff daily; monthly counts within floor/ceiling limits; no consecutive-day assignment. | Hard |
| Managerial coverage | At least 1 section chief or higher overall and at the main office; 2 in total preferred. | Hard / Soft |
| Individual and fixed schedules | Employment-specific days off, precedence of fixed days off and events, and part-time work restrictions. | Hard within defined scope |
| Weeks and horizon | Weekly days off evaluated over full weeks; explicit horizon-edge and individual eligibility exceptions. | Hard in Stage 1; relaxed in Stage 2 |
| Consecutive conditions | Limits on consecutive effective-work and roster-state days, and consecutive ordinary days off. | Hard |
| Fairness | Hard monthly branch-duty limits and full-horizon holiday/weekend day-off ranges (≤1); additional monthly range/deviation and repeat-duty penalties. | Hard / Soft |

Note: Weekly exceptions and Stage 2 relaxation are defined in Sections 3.6 and 3.11; fairness includes hard and soft components.

## 2.3 Manual roster preparation and its limitations

Previously, staff entered fixed schedules into a spreadsheet and manually adjusted and checked headcounts, managerial coverage, branch-office assignments, days off and consecutive working days. A change to one assignment could affect another requirement, necessitating repeated adjustment and checking.

We estimated that preparing a monthly roster required 12–18 h of work, spread over 2–3 calendar days including staff review (Section 4.5). Checking month-boundary weeks and consecutive working days was particularly burdensome. Adjustments depended on the experience of a staff member familiar with the work rules, raising concerns about continuity during absence or reassignment.

The proposed system represents these requirements explicitly, models the entire planning horizon and generates a draft roster with information for checking.

## 3 METHODS

This section formulates the department's work rules in terms of the constraints and objective used in the implementation. It distinguishes optimisation over the entire horizon, exceptions to weekly day-off rules, conditional two-stage solving and the generation of checking information.

### 3.1 Sets and decision variables

Let $N(S)$ denote the cardinality of a set $S$. Let $I$ be the staff set and $D$ the chronologically ordered set of planning dates, with $n = N(I)$ and $T = N(D)$. The sets of regular, contract (shokutaku) and part-time staff are denoted by $I_R, I_C, I_P$, respectively, and the set of section chiefs or higher by $I_M$. Let $D_m$ contain the planning dates in month $m$. For staff member i, $a_{i,d} = 1$ if date $d$ falls within that person's active period, and zero otherwise. An active period defines eligibility for assignment, for example by employment dates; it does not mean that the person actually works on that day.

Let $D_S$ be the dates registered in the input file for public holidays and office-closed days. Saturdays, Sundays and dates in $D_S$ form $D_{hd}$; the remaining dates form $D_{wd}$. These input sets determine effective-headcount requirements and day-off accounting. Section 3.12 explains the relationship between the official holiday calendar and the input file.

Three binary variables are defined for each staff–day pair:

$w_{i,d}$: roster state for work and schedule management, equal to one on ordinary working days, non-working event days and fixed holiday days off, and zero on ordinary weekly days off and other non-working days such as public holidays.

$e_{i,d}$: one when the staff member counts towards the department's effective headcount, and zero otherwise.

$k_{i,d}$: one when the staff member is assigned to the branch office, and zero otherwise.

The variable $w_{i,d}$ is neither an effective-work indicator nor a direct representation of the roster's day-duty symbol. For example, paid leave, training and business trips have $w_{i,d} = 1, e_{i,d} = 0$. Fixed holiday days off have the same values but are displayed as holiday-off entries. Weekly days off are counted by summing $1 - w_{i,d}$ over the eligible dates defined in Section 3.6. This separates schedule management from effective staffing.

### 3.2 Consistency between variables

Branch-office duty requires both a working roster state and effective work:

$$k_{i,d} \le w_{i,d}, \qquad k_{i,d} \le e_{i,d} \qquad (i \in I, d \in D). \qquad (1)$$

Outside the active period, $w_{i,d} = e_{i,d} = 0$, and hence $k_{i,d} = 0$ by Equation (1). Within the active period, an ordinary day without a fixed schedule, registered event or group meeting satisfies

$$e_{i,d} = w_{i,d}. \qquad (2)$$

Other dates use the fixed values by event type in Section 3.8. These rules imply $e_{i,d} \le w_{i,d}$ for every staff–day pair. Equation (2) also prevents the model from assigning an arbitrary working roster state without counting the person as effective staff on ordinary days.

### 3.3 Effective headcount

Ordinary weekdays require a minimum headcount, whereas weekends, public holidays and office-closed days require an exact headcount. Both requirements are six staff in this implementation:

$$\sum_{i \in I} e_{i,d} \ge 6 \quad (d \in D_{wd}), \qquad \sum_{i \in I} e_{i,d} = 6 \quad (d \in D_{hd}). \qquad (3)$$

The six include staff assigned to the branch office. Non-working events such as leave, training and business trips are excluded from effective staffing through the rules in Section 3.8.

### 3.4 Managerial coverage

At least one section chief or higher must be available each day, with at least one at the main office:

$$\sum_{i \in I_M} e_{i,d} \geq 1, \qquad \sum_{i \in I_M} \left(e_{i,d} - k_{i,d}\right) \geq 1 \quad (d \in D). \qquad (4)$$

The implementation adds these constraints when the input includes at least one section chief or higher. The evaluation instance contains eight such staff. A target of two managers is distinguished from the mandatory minimum of one; its shortfall enters the objective in Section 3.10.

### 3.5 Branch office assignments

Exactly two staff are assigned to the branch office each day:

$$\sum_{i \in I} k_{i,d} = 2 \quad (d \in D). \qquad (5)$$

If a date has an entry in a staff member's individual event register, including a fixed day off, $k_{i,d} = 0$. Group-meeting settings are separate from that register, so a meeting day alone does not prohibit branch duty. A staff member cannot be assigned to the branch on two consecutive dates within their active period:

$$k_{i,d} + k_{i,d+1} \leq 1. \qquad (6)$$

For month $m$, let $I_m$ be the staff with at least one active date and $n_m$ their number. Monthly branch-duty counts $K_{i,m} = \sum_{d \in D_m} k_{i,d}$ satisfy

$$\text{floor}\left(\frac{2N(D_m)}{n_m}\right) \leq K_{i,m} \leq \text{ceil}\left(\frac{2N(D_m)}{n_m}\right) \quad (i \in I_m). \qquad (7)$$

Here, floor and ceil denote the floor and ceiling functions. Equation (7) is a hard balancing constraint across eligible staff in each month, irrespective of employment category. It does not allocate duty in proportion to each person's number of active dates. Full-horizon imbalance and repeat assignments after two or three days are penalised separately as soft constraints.

### 3.6 Weekly days off and exceptions

For staff member $i$, let $c_i$ be the standard weekly working-day count and $r_i = 7 - c_i$ the standard weekly day-off count. The values are $(c_i, r_i) = (5,2)$ for regular staff, $(4,3)$ for contract staff and $(3,4)$ for part-time staff. Weeks $t$ are partitioned according to the specified starting weekday. Let $\mathcal{T}$ be the weeks whose seven dates all lie within the planning horizon.

For staff member $i$ and week $t$, let $A_{i,t}$ be the dates eligible for weekly day-off accounting. A fixed weekly day off within the active period is included in $A_{i,t}$ even when it coincides with a holiday, office closure or another event. Other dates are included only if they are within the active period and are neither holidays or office-closed days nor group meetings, non-working events or other registered events. Giving fixed weekly days off precedence is essential to avoid adding too many days off in the same week.

Let $f_{i,t}$ be the number of fixed weekly days off in $A_{i,t}$. Required and actual day-off counts are defined as

$$r_{i,t} = \max(r_i, f_{i,t}), \qquad q_{i,t} = \sum_{d \in A_{i,t}} (1 - w_{i,d}). \qquad (8)$$

In Stage 1, the weekly equality applies only to the set $\mathcal{K}_1$ of staff–week pairs in full weeks satisfying $N(A_{i,t}) \geq r_i$:

$$q_{i,t} = r_{i,t} \quad ((i,t) \in \mathcal{K}_1). \qquad (9)$$

There are two exceptions. First, partial weeks at the horizon edges are excluded from weekly day-off evaluation in both stages. Second, even in a full week, Stage 1 omits the equality for staff–week pairs with fewer candidate dates than the standard $r_i$. This second exception is based on candidate counts; it does not establish compatibility with staffing requirements. A Stage 1 solution therefore does not guarantee the absence of shortages for exempted staff–week pairs. After solving, shortages are recounted for every full week and reported in a checking sheet.

Stage 2 evaluates shortages over $I \times \mathcal{T}$, including staff–week pairs covered by the second exception. Section 3.11 defines the scope and relaxation used in each stage.

For part-time staff, $w_{i,d} = e_{i,d} = 0$ is imposed separately on weekends, public holidays and office-closed days. Let $E_{i,t}$ contain the active dates within a week excluding weekends, holidays, office-closed days and non-working events. The following applies, including in partial weeks:

$$\sum_{d \in E_{i,t}} e_{i,d} = \min\left(3, N\left(E_{i,t}\right)\right) \quad (i \in I_P). \qquad (10)$$

Thus, when fewer than three eligible dates are available, their number is the upper bound. This part-time working-day condition is retained in Stage 2.

### 3.7 Consecutive work and days off

Regular staff may perform effective work on at most five dates in any six consecutive dates. For all staff, both effective-work days and the sum of roster states $w$ are limited to six in any seven consecutive dates. Each window is evaluated only when all its dates lie within the planning horizon:

$$\sum_{j=0}^{5} e_{i,d+j} \leq 5 \quad (i \in I_R), \qquad (11)$$

$$\sum_{j=0}^{6} e_{i,d+j} \leq 6, \qquad \sum_{j=0}^{6} w_{i,d+j} \leq 6 \quad (i \in I). \qquad (12)$$

The condition on $w$ limits consecutive roster states, including events that do not count as effective work. A fixed holiday day off has $w = 1$ but is displayed as a holiday-off entry. This condition therefore differs from counting consecutive work solely from the displayed roster symbols.

For regular staff, if three consecutive dates all lie within the active period and none is a public holiday or office-closed day,

$$w_{i,d} + w_{i,d+1} + w_{i,d+2} \geq 1 \quad (i \in I_R). \qquad (13)$$

This prohibits three consecutive ordinary days off. Equation (13) is not applied when the window contains a public holiday or office-closed day.

### 3.8 Fixed schedules and events

Daily states follow Table 3 from top to bottom. ‘Other’ means that no higher-priority rule applies. Dates outside the active period take precedence over all other conditions.

***Table 3. Daily state assignments and rule precedence.***

| Priority and condition | $w_{i,d}$ | $e_{i,d}$ | Treatment |
|---|---|---|---|
| 1. Outside the active period | 0 | 0 | Not assigned; blank roster entry. |
| 2. Fixed weekly day off | 0 | 0 | Counted as a weekly day off; no branch duty. |
| 3. Fixed holiday off | 1 | 0 | Not counted as a weekly day off; shown as a holiday-off entry. |
| 4. Group-meeting day | 1 | 1 | Overrides non-working events, except for the higher-priority fixed days off. |
| 5. Non-working event | 1 | 0 | Annual, summer or special leave; training; business trip; morning/afternoon leave. |
| 6. Other registered event | 1 | 1 | Committee meetings, works, etc.; no branch duty. |
| 7. Ordinary day without a fixed event | Free | $w_{i,d}$ | Work and days off are decided by optimisation. |

Non-working events and fixed holiday days off are distinguished from weekly days off in roster management. Half-day leave is treated as a full-day non-working event for effective-headcount accounting in this implementation. If a non-working event is registered on a meeting day, the meeting rule takes precedence except when a fixed weekly day off or fixed holiday day off applies.

A fixed weekly day off on a public holiday or office-closed day still counts towards weekly days off, although the roster's first row displays a holiday-off entry. Consequently, the checking sheet's weekly count may differ from a count of displayed weekly-day-off symbols under some inputs. Weekly days off are evaluated using Equation (8).

### 3.9 Optimisation over the full horizon and month boundaries

A single model is constructed over the date set $D$ for the entire horizon. A week spanning two months is treated as seven consecutive dates, with weekly day-off rules applied to staff–week pairs meeting the conditions in Section 3.6. Consecutive-work rules and branch-duty intervals are also evaluated across month boundaries.

If months are solved independently without carrying information across boundaries, a boundary week becomes a partial week in each model and is excluded from the weekly equality. The combined roster then lacks a constraint on the total days off over that full week. Figure 1 illustrates the difference.

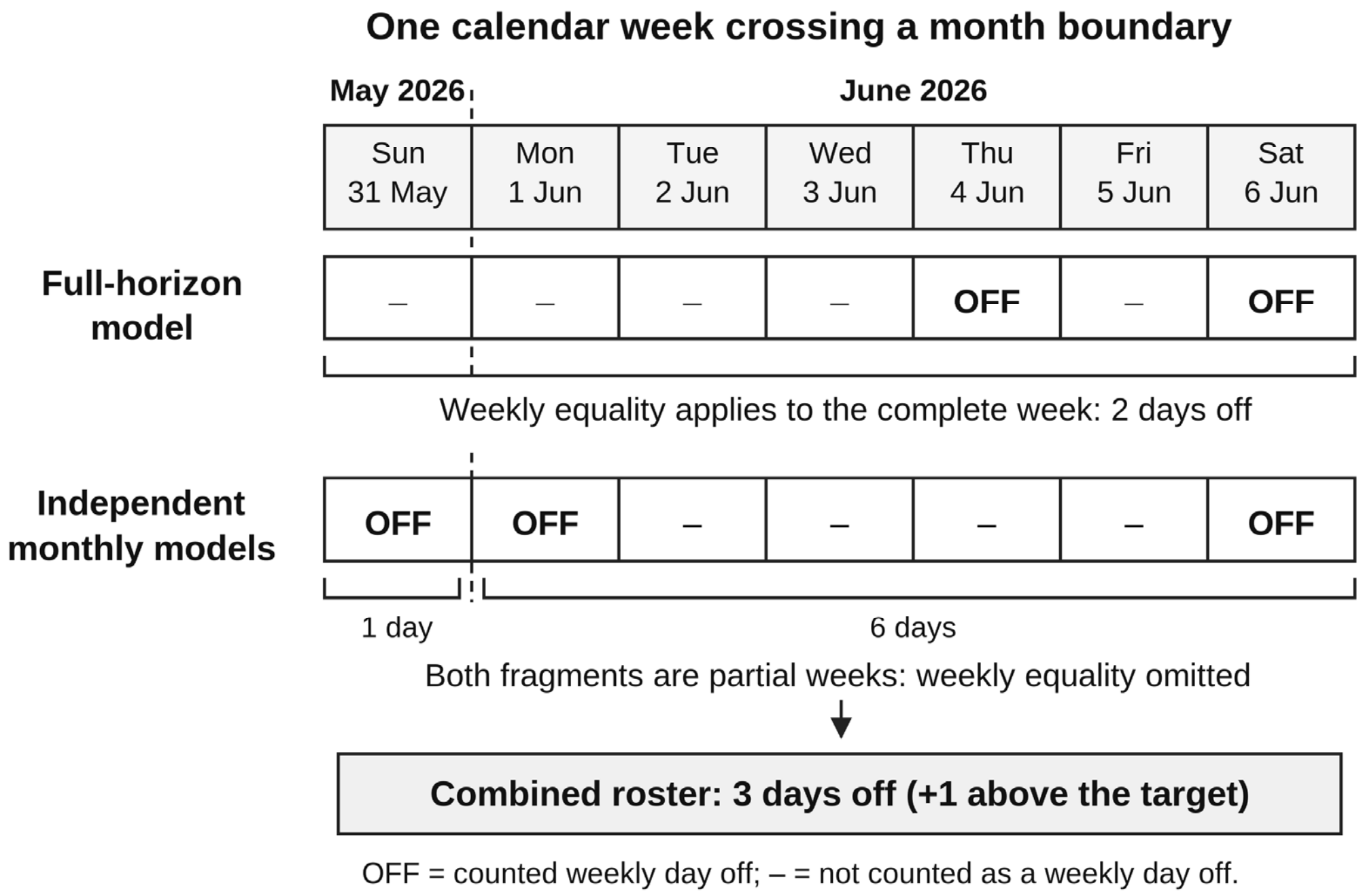

***Figure 1. Weekly day-off accounting across a month boundary: example from 31 May–6 June 2026.***
Allocations correspond to the example in Section 4.7. The full-horizon equality applies within the scope of Eq. (9). In independent monthly models, each fragment is a partial week and the weekly equality is omitted; other constraints and objectives remain active.

Decomposition with assignment history can also evaluate period boundaries (Kingston, 2021). This study uses a full-horizon formulation and evaluates weeks and consecutive-day windows wholly contained in that horizon. No assignments outside the horizon are supplied.

## 3.10 Hard and soft constraints and the objective function

The distinction between hard and soft constraints reflects mandatory and desirable conditions identified with the department, rather than whether a condition can technically be written as an equality or inequality. For example, one manager is mandatory and two are a target. A desirable condition could be made mandatory by modifying the model, for example by fixing its shortfall to zero, but compatibility with other requirements would need to be reassessed. This is a possible model modification, not a user-interface option for freely switching constraint types.

Balancing uses both hard and soft constraints. In addition to monthly branch-duty counts in Equation (7), the combined set of regular and contract staff $I_{RC} = I_R \cup I_C$ is subject to the following conditions on full-horizon non-working counts on holidays and office-closed days $H_i$ and ordinary weekend day-off counts $P_i$:

$$\max_{i\in I_{RC}} H_i - \min_{i\in I_{RC}} H_i \leq 1, \qquad \max_{i\in I_{RC}} P_i - \min_{i\in I_{RC}} P_i \leq 1. \qquad (14)$$

The count $H_i$ sums $1 - w_{i,d}$ over active dates in $D_S$. The count $P_i$ sums $1 - w_{i,d}$ over active Saturdays and Sundays excluding public holidays, office-closed days and non-working events. Equation (14) includes staff with no eligible dates as having zero counts. Fixed holiday days off have $w = 1$ and are therefore excluded from $H_i$, despite their holiday-off display.

The soft-constraint measures are defined next. For a vector of counts $x = (x_i)_{i\in G}$, define the range and an integer-scaled absolute-deviation measure as

$$R_G(x) = \max_{i\in G} x_i - \min_{i\in G} x_i. \qquad (15a)$$

$$A_G(x) = \sum_{i\in G} \phi\left(N(G)x_i - \sum_{j\in G} x_j\right). \qquad (15b)$$

where $\phi(x) = \max(x, -x)$ denotes absolute value. The measure $A_G$ is $N(G)$ times the sum of absolute deviations from the mean. Balancing therefore uses absolute deviations for some terms as well as maximum–minimum differences.

The monthly non-working count on holidays and office-closed days $h_{i,m}$ sums $1 - w_{i,d}$ over dates in $D_S \cap D_m$ that lie within the active period and are not non-working event days. Counts are compared separately for regular and contract staff combined and for part-time staff. Staff with no eligible dates in that month are excluded. Denote a comparison group by $G$.

Let $p_{i,m,b}$ be monthly day-off counts for $b = \text{Sat}, \text{Sun}, \text{Weekend}$. Eligible dates are active dates of the relevant weekday, excluding holidays, office-closed days and non-working events; their number is $C_{i,m,b}$. Comparisons combine regular and contract staff and include only those with $C_{i,m,b} > 0$. Each person's target interval is

$$L_{i,m,b} = \text{floor}(C_{i,m,b}/2), \qquad U_{i,m,b} = \min(C_{i,m,b}, L_{i,m,b} + 1) \qquad (16)$$

with deviations

$$u^L_{i,m,b} = \max(0, L_{i,m,b} - p_{i,m,b}), \qquad u^U_{i,m,b} = \max(0, p_{i,m,b} - U_{i,m,b}) \qquad (17)$$

Because the upper bound is based on the lower bound plus one, the interval for an even number of eligible dates differs from fixing both bounds at the same floor/ceiling value.

Deviation from contractual working-day counts is calculated for every staff–week pair, including partial weeks. Let $B_{i,t}$ be the active dates within a week and $\tau_{i,t} = \mathrm{round}\big(c_i N(B_{i,t})/7\big)$ the prorated target. Rounding is half up, and the target is bounded between zero and $N(B_{i,t})$. The deviation is $v_{i,t} = \phi\big(\sum_{d\in B_{i,t}} e_{i,d} - \tau_{i,t}\big)$.

Let $K_i = \sum_{d\in D} k_{i,d}$ be full-horizon branch-duty counts. As two staff are assigned daily, the total count is $2T$, giving the absolute-deviation term $\sum_{i\in I} \phi\,(nK_i - 2T)$. Repeat assignments are represented by $z_{i,d}^{(\ell)} = k_{i,d}k_{i,d+\ell}$ for $\ell = 2{,}3$, where both dates are active. The implementation expresses each logical conjunction using three linear inequalities.

The objective $F$ is the weighted sum of the penalty measures in Table 4. The weights are the implemented values, and all measures use binary or integer variables and are evaluated as integers.

***Table 4. Active objective terms and implemented weights.***

| Objective term | Scope and penalty measure | Weight |
|---|---|---|
| Shortfall from 2 managers | All dates; sum of $\max\big(0{,}2 - \sum_{i\in I_M} e_{i,d}\big)$ | 1 |
| Contract-work deviation | All staff and weeks; sum of $v_{i,t}$. Partial weeks use prorated targets. | 500 |
| Monthly holiday-off deviation | By month and comparison group: $A_G\big(h_{\cdot,m}\big)$ | 20,000 |
| Excess monthly holiday-off range | By month and comparison group: $\max\big(0, R_G(h_{\cdot,m}) - 1\big)$ | 400,000 |
| Monthly Saturday day-off range | Range $R$ among eligible regular and contract staff. | 120,000 |
| Monthly Sunday day-off range | Same scope as the Saturday range. | 120,000 |
| Monthly total weekend day-off range | Same scope as the Saturday range. | 120,000 |
| Individual monthly weekend targets | For Saturday, Sunday and both combined: $u^L + u^U$ | 180,000 |
| Full-horizon branch-duty range | All staff: $R_I(K)$ | 1,000 |
| Full-horizon branch-duty deviation | All staff: $\sum_i \phi\,(nK_i - 2T)$ | 1 |
| Branch duty repeated after 2 days | Sum of $z_{i,d}^{(2)}$ | 30 |
| Branch duty repeated after 3 days | Sum of $z_{i,d}^{(3)}$ | 15 |

| Objective term | Scope and penalty measure | Weight |
| --- | --- | --- |
| Roster state on holidays / closed days | Sum of $w_{i,d}$ over active staff–day pairs. | 50 |

For monthly count comparisons, range and mean-deviation terms are omitted when the comparison group has at most one member. Individual monthly terms are also omitted if there are no eligible dates. The final holiday/closure term evaluates $w$, not $e$; it therefore does not directly minimise effective staffing. Full-horizon branch-duty comparisons include all staff without prorating for employment category or active-period length. Monthly balancing of total days off and penalties for violating requested days off are not active objective terms in this implementation.

The objective weights were established through an iterative calibration process involving repeated computational experiments and consultation with the staff member responsible for roster preparation. Because no established monetary cost was available for the different scheduling preferences, the weights were used to represent their relative operational importance rather than economic values. Fairness-related measures and monthly balance criteria were assigned substantially larger weights than lower-priority preferences, such as branch-duty spacing and managerial-preference targets, to ensure that improvements in secondary objectives would not be obtained at the expense of conditions considered more important in practice. Candidate weight sets were tested using operational data, and the final values were retained when the resulting rosters were judged by operational staff to provide an acceptable balance between feasibility, fairness and practical usability. The selected weight set also produced stable solution behaviour across the test instances examined during system development.

### 3.11 Conditional two-stage solving and solver statuses

Stage 1 minimises $F$ subject to the normal model's hard constraints, including the exception rules in Section 3.6. A solution is output when the status is OPTIMAL or FEASIBLE. OPTIMAL indicates proven optimality; FEASIBLE indicates that a feasible solution has been found without proof of optimality (Perron & Didier, 2025).

Stage 2 is a second solver call that attempts to obtain a reference roster identifying adjustments when the normal model yields no feasible solution. It shifts the system from ordinary roster generation to support for revising input conditions.

The implemented trigger is any Stage 1 status other than OPTIMAL or FEASIBLE. It therefore includes UNKNOWN, for example after reaching a time limit, as well as INFEASIBLE. INFEASIBLE proves that no feasible solution exists; UNKNOWN does not determine whether one exists. The current branch also routes MODEL_INVALID to Stage 2. This is a model-validation error and must not be interpreted as infeasibility of the operational requirements (Perron & Didier, 2025).

Stage 2 considers $\mathcal{K}_2 = I \times \mathcal{T}$ over full weeks and introduces integer shortage variables $s_{i,t}$. This includes staff–week pairs excluded from the Stage 1 equality because of insufficient candidate dates:

$$q_{i,t} \le r_{i,t}, \qquad q_{i,t} + s_{i,t} \ge r_{i,t}, \qquad 0 \le s_{i,t} \le r_{i,t} \quad \big((i,t) \in \mathcal{K}_2\big). \qquad (18)$$

All other hard constraints are retained, and the objective is

$$\min\left(M \sum_{(i,t)\in\mathcal{K}_2} s_{i,t} + F\right), \qquad M = 1{,}000{,}000{,}000. \qquad (19)$$

To penalise day-off shortages strongly, the total shortage receives weight $M = 10^9$ in Equation (19). OPTIMAL indicates optimality for this weighted objective; a FEASIBLE status yields a feasible solution without that proof.

If Stage 2 finds a zero-shortage solution, Equation (18) implies $q_{i,t} = r_{i,t}$ for every relevant pair, and the Stage 1 constraints are also satisfied. A zero-shortage relaxed solution therefore establishes the existence of a feasible solution to the normal model. Conversely, Stage 1 success alone does not rule out shortages for pairs exempted from Equation (9).

Stage 2 output is a reference roster for reviewing inputs and schedules. If the retained headcounts, fixed days off, managerial coverage, branch assignments or other hard conditions are incompatible, Stage 2 can also fail. In that case, the program terminates without generating a new roster file.

### 3.12 Implementation and checking information

The system is implemented in C# and solved using CP-SAT in Google OR-Tools (Perron et al., 2023; Perron & Didier, 2025). ClosedXML produces the Excel roster. After installation of the required software, input from comma-separated value (CSV) files, solving and output can be completed on a single personal computer without an external server connection.

Because the main decisions concern daily assignments and headcounts, the formulation primarily uses binary variables and integer linear constraints. Maximum, minimum and absolute-value expressions use CP-SAT’s AddMaxEquality, AddMinEquality and AddAbsEquality constraints. This differs from formulations primarily concerned with overlapping intervals or the starting times of variable-duration tasks.

The main inputs are staff records including employment category, grade and active period; holidays and office-closed dates; and individual fixed schedules and events. The horizon, starting weekday and ordinal weekday for meetings are specified at runtime. The program also reads a requested-day-off CSV file, but those requests are not used by this version’s constraints or objective.

Headcount requirements, holiday classifications and calendar rules for meeting transfers depend on the input CSV files. The public-holiday set generated within the program is used mainly for roster colouring and does not fill missing CSV dates. Consistency between the public calendar and the CSV is an input-data management responsibility.

Meetings default to the third Thursday of each month; ‘third’ denotes that weekday’s third occurrence in the month. A meeting on a non-business day is transferred before solving to the first subsequent business day within the horizon. If no such date exists within the horizon, it is recorded as unset. The transfer date in the output is the setting already applied to the model.

After a solution is found in either stage, the system outputs a ‘Weekly two-days-off check’ sheet. Despite that label, it uses employment-specific requirements to count shortages for every staff–week pair in full weeks. It shows required and actual days off, shortage, fixed weekly days off, candidate dates and remarks, or indicates that no shortage exists. Remarks are standard categories based on whether candidate dates are insufficient.

When meetings are transferred, a ‘Group-meeting transfers’ sheet records the original date, the applied date or an unset status, and the calendar-based reason. Figure 2 shows the input, solving and output process.

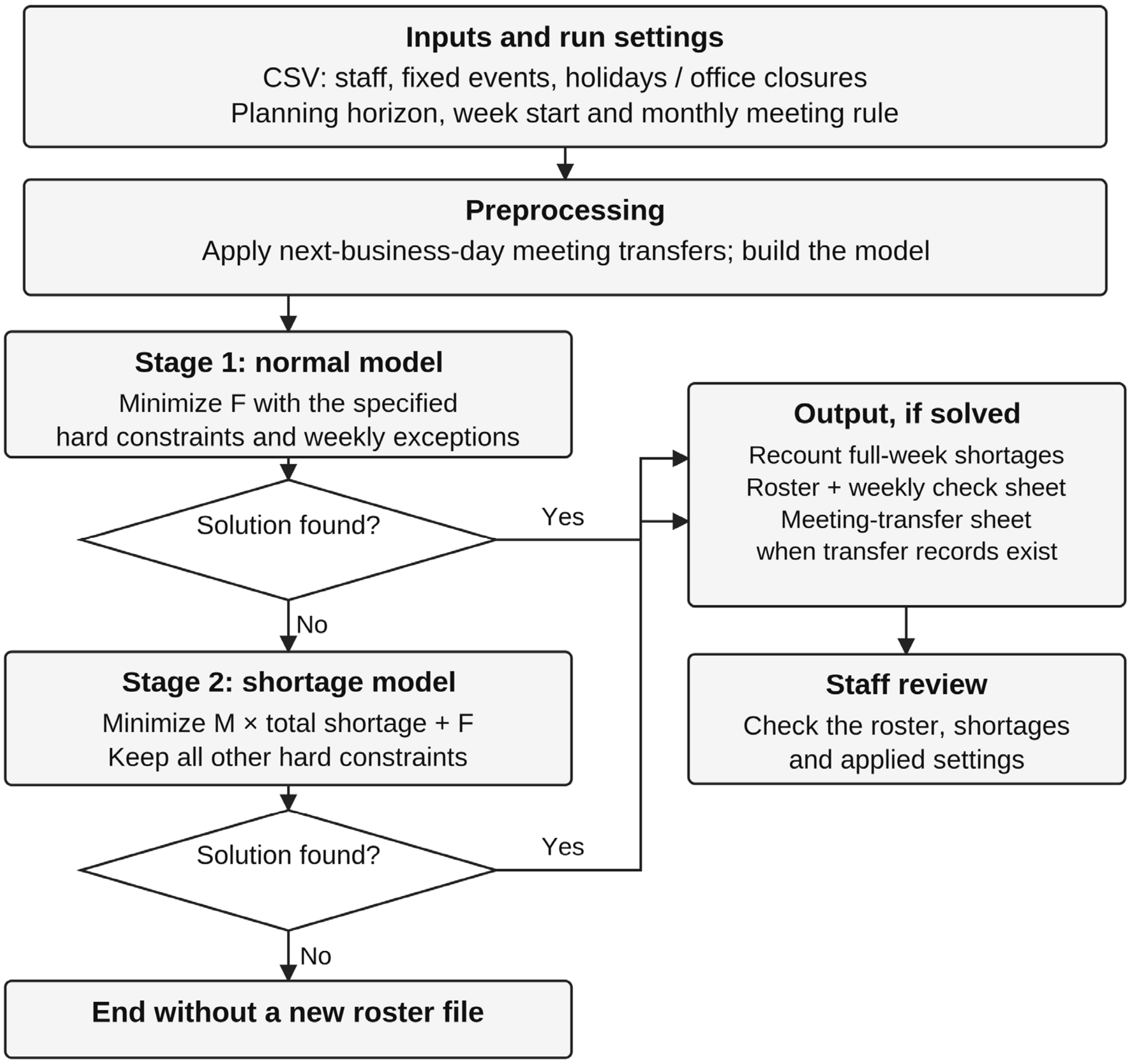

***Figure 2. Input, conditional two-stage solving and output of the implemented system.***

"Solution found" denotes OPTIMAL or FEASIBLE; other statuses follow the next branch (Section 3.11). Full-week shortages are recounted after either successful stage. Figure labels are English translations of the deployed Japanese interface.

Simple pre-solve checks display warnings about daily and weekly input conditions in the standard output. Feasibility and optimality are determined from the solver status for the relevant model.

# 4 DEPLOYMENT AND EVALUATION RESULTS

## 4.1 Operational scope and evaluation design

This section reports deployment in routine operations and the estimated change in roster preparation time (Section 4.5). Computational experiments assess system behaviour under the baseline (Section 4.3), a modified meeting-day scenario (Section 4.6), and alternative partitions of the planning horizon (Section 4.7). Inputs comprised staff, calendar and fixed-schedule CSV files.

The computational horizon was 1 May–31 December 2026, totalling 245 days. The system accepts specified start and end dates, and these dates were selected for the experiment. The roster covered 14 staff in the Water Purification Office: eight section chiefs or higher and six other staff, comprising 13 regular and one contract employee. It included no part-time staff and contained 10 fixed schedule entries. Weeks began on Sundays, and baseline meetings were set to the third Thursday of each month. The official calendar contains 12 public holidays or other statutory holidays in this period (Cabinet Office, 2026). The holiday/closure set used in solving nevertheless depends on the input CSV, as described in Section 3.12.

## 4.2 Model size

Table 5 gives the instance size and computing environment. There were 20,968 variables including auxiliary variables, of which 17,173 were binary. The listed constraint counts are linear constraints; the model also uses maximum, minimum and absolute-value constraints.

***Table 5. Planning instance, model size and computing environment.***

| Item | Value / setting |
|---|---|
| Number of staff | 14 |
| Section chief or higher | 8 |
| Regular employees | 13 |
| Contract (shokutaku) employees | 1 |
| Planning horizon | 1 May - 31 Dec 2026 |
| Number of days | 245 |
| Fixed individual schedules | 10 |
| Decision variables | 20,968 |
| Boolean (0-1) variables | 17,173 |
| Linear constraints (total) | 51,897 |

| Item | Value / setting |
|---|---|
| Two-variable linear constraints | 27,922 |
| Three-variable linear constraints | 9,855 |
| Multi-variable linear constraints | 12,921 |
| Single-variable linear constraints | 1,199 |
| Processor (CPU) | 12th Gen Intel Core i7-1255U |
| Operating system | Windows 11 |
| Implementation language | C# |
| Solver | CP-SAT (Google OR-Tools v9.14.6206) |
| Parallel workers | 8 |
| Time limit per solver call | 600 s |

### 4.3 Baseline solution and performance

With eight parallel workers and a 600 s limit per solver call, the baseline model for 14 staff and 245 days reached OPTIMAL in approximately 177 s, with matching objective and lower-bound values. The solution minimised the objective in Table 4 and satisfied the hard constraints included in the model.

The baseline finished in Stage 1. Post-solve recounting confirmed zero total day-off shortage over full weeks, within the scope defined in Section 3.6.

### 4.4 Horizon edges and partial weeks

With weeks starting on Sunday, the first week is 26 April–2 May 2026, of which only 1–2 May lies within the horizon. The final week is 27 December 2026–2 January 2027, of which 27–31 December lies within the horizon. Both are excluded from weekly day-off evaluation in both stages under Section 3.6. Full weeks spanning internal month boundaries remain subject to the specified weekly rules.

### 4.5 Operational deployment and reported effects

Operational use began at the Water Purification Office of Suita City Waterworks Bureau in May 2026. The system generates an Excel draft roster and checking information from CSV inputs, followed by staff review and adjustment as part of the normal preparation workflow.

The reported times are operational estimates provided by the staff member responsible for roster preparation, who had experience preparing monthly rosters both before and after deployment of the system (Figure 3). The estimates refer to the total working time required to prepare a monthly roster, including confirmation of leave requests and other fixed schedules, roster generation, staff review, revisions and final checking. They were derived from operational experience accumulated during routine use of the system rather than from a formal time-and-motion study. The largest reduction was associated with the assignment of working days, days off and branch-office duties, which had previously required repeated manual adjustment and verification.

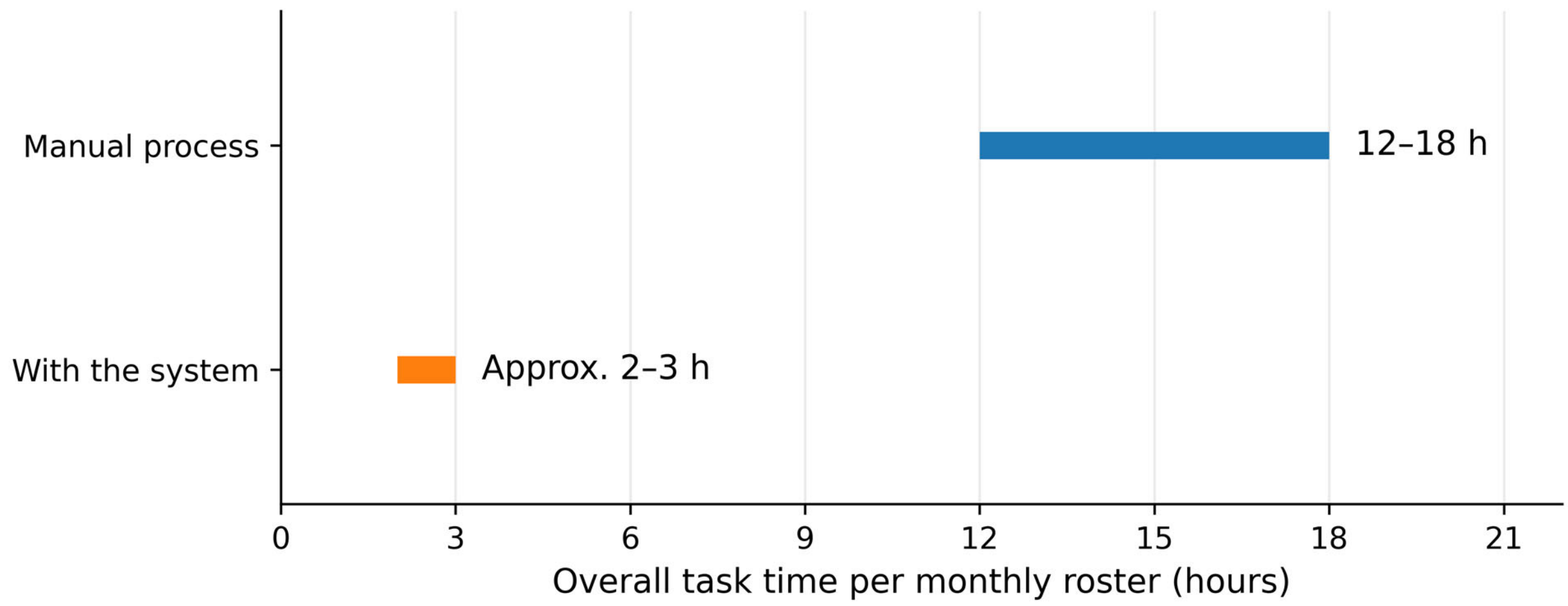

***Figure 3. Estimated overall task time before and after operational deployment.***

Lines show estimated ranges of overall task time, rather than solver time. The manual process extended over 2–3 calendar days.

Table 6 gives the estimated breakdown of working time. Assigning working days and days off decreased from 8–10 h to at most 1 h. Draft presentation, staff review, revisions and final checking still required 1 h after deployment.

*Table 6. Estimated working time per monthly roster before and after deployment.*

| Work stage | Before deployment (manual) | After deployment (system-assisted) |
|---|---|---|
| Confirmation and entry of fixed schedules (leave requests, business trips, training, works, committees) | 1–2 h | up to 1 h |
| Assignment of working days and days off | 8–10 h | up to 1 h |
| Selection of branch-office duty staff (two per day) | 1–2 h | Not required |
| Checking the entered content | 1–2 h | Not required |
| Draft presentation, staff review, revisions, and final check | 1–2 h | 1 h |
| Total | 12–18 h (over 2–3 calendar days) | approx. 2–3 h |

Note: Times are operational estimates (Section 4.5). The manual total extends over 2–3 calendar days, including waiting for staff to review the draft. The total compares the overall task time, not the solver computation time. “Not required” refers only to the indicated stages, not to all checking or review activities.

## 4.6 Comparison with and without the second stage

Meetings were changed from the third Thursday to the third Monday of every month. We compared ending after the normal model with proceeding to Stage 2 when that model yielded no feasible solution.

Under this setting, Stage 1 returned INFEASIBLE in approximately 3.7 s.

Stage 2 reached OPTIMAL for the weighted objective in Equation (19) in approximately 244 s. The reference roster retained the other hard constraints. During 20–26 September 2026, five staff had one weekly day off against a requirement of two, giving five staff-days of shortage in total.

The official calendar lists holidays on 21, 22 and 23 September that week (Cabinet Office, 2026). If these are also registered in the CSV, the 21 September meeting is transferred to 24 September, affecting candidate dates for weekly days off. Actual candidate counts depend on the calendar input and each person’s fixed schedules.

Figure 4 compares statuses and shortages under the two meeting settings. The approximately 244 s is the Stage 2 solving time alone and excludes the approximately 3.7 s spent in Stage 1.

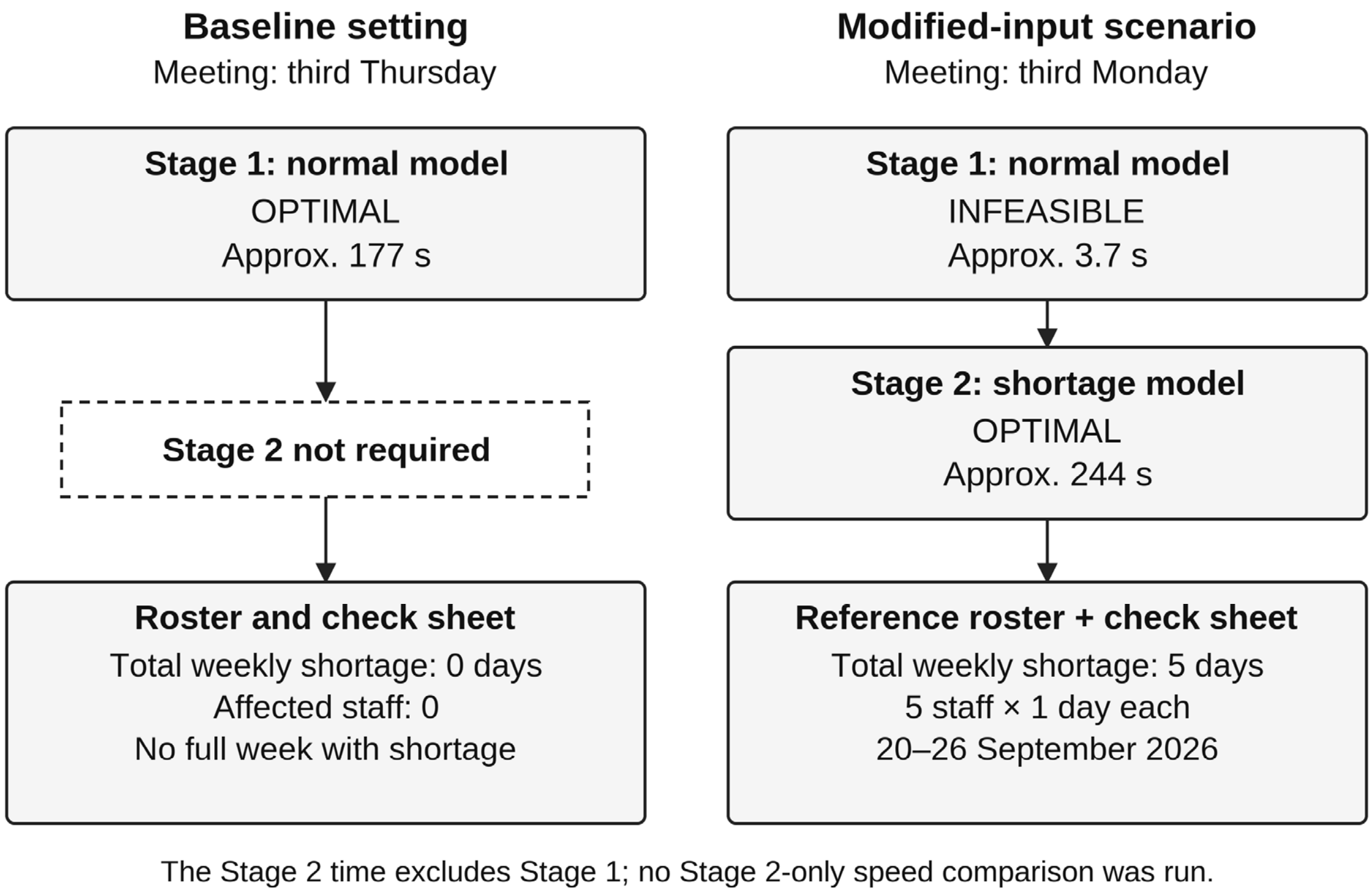

***Figure 4. Two-stage results under the baseline and modified group-meeting rules.***

In the modified-input scenario, five staff each received 1 day off against a requirement of 2 in the week of 20–26 September 2026. Stage 2 time excludes Stage 1.

### 4.7 Full horizon and independent monthly solving

We compared solving the full horizon with solving May–December 2026 independently by month and concatenating the results. Operational inputs and the third-Thursday meeting setting were shared. All eight monthly solves reached optimality in Stage 1. The independent approach carried neither assignment history nor full-horizon totals across months, and we evaluated weekly days off at the resulting boundaries.

Among six full Sunday-start weeks spanning month boundaries, independent monthly solving gave regular staff three rather than two weekly days off in two weeks. There were 16 staff–week cases and 16 excess staff-days in total. Because each month's boundary fragment is excluded from the weekly equality, the combined week's day-off count is not constrained.

For 31 May–6 June 2026, the full-horizon example assigns two days off, on 4 and 6 June. Independent monthly solving assigns three, on 31 May and 1 and 6 June (Figure 1). The full-horizon model applies the weekly equality to the same seven dates.

The staff count in Figure 5 is the number whose weekly days off under independent monthly solving exceed those under full-horizon solving. Each affected staff–week has one excess day; zero indicates no excess under this measure.

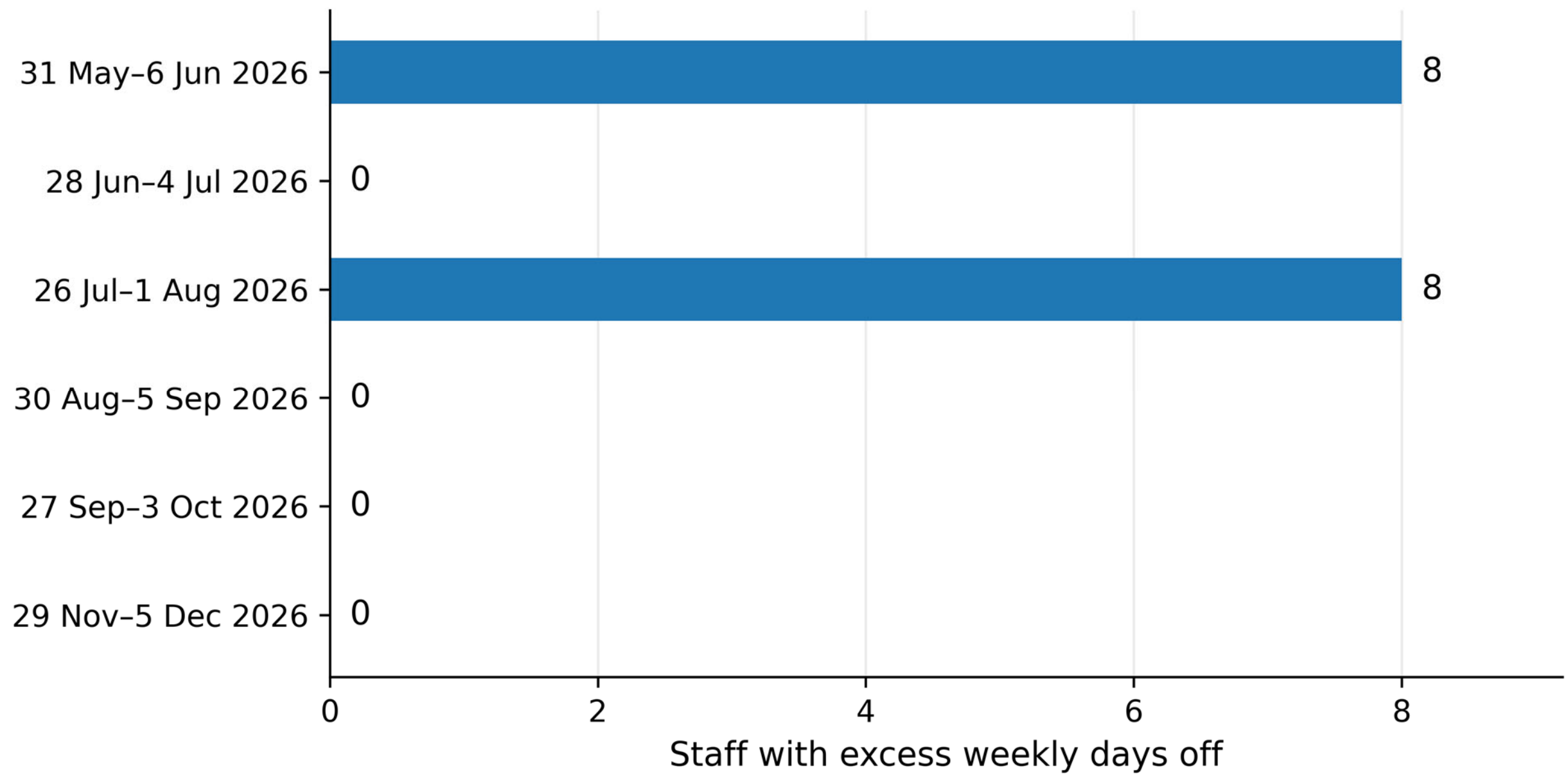

***Figure 5. Staff with excess weekly days off under independent monthly solving.***
Affected staff received 3 days off versus 2 under full-horizon solving. The total is 16 staff–week cases and 16 extra staff-days, not 16 distinct people. Zero indicates no excess under this measure.

### 4.8 Summary of results

Operational estimates indicate reduced preparation time, particularly for assigning work and days off. The baseline computational experiment obtained a Stage 1 roster, while Stage 2 generated a reference roster identifying shortages after the meeting-day change. The independent-month comparison demonstrated the need for accounting across month boundaries.

## 5 DISCUSSION

### 5.1 System design and performance in support of operational use

The system combines CSV input, solving on a single computer and Excel output to fit the department's preparation and checking procedures. Its practical contribution includes both modelling the work rules and integrating spreadsheet-compatible input and output with information that staff can inspect. Solving the baseline within a few minutes provides evidence of the computational performance supporting this design. Table 5 specifies the environment and stopping conditions.

## 5.2 Practical role of conditional second stage solving

In the modified meeting-day scenario, Stage 2 returned a reference roster with specific adjustment needs after the normal model proved infeasible. Compared with terminating at that point, reporting shortages by staff member and week gives the responsible person information for reviewing the inputs.

The system ordinarily generates a roster under the specified mandatory conditions and attempts reference-roster generation if no solution is obtained. This separation places routine roster preparation and support for revising requirements within one workflow.

## 5.3 Full horizon optimisation and month boundaries

The excess weekly days off under independent monthly solving without history illustrate a boundary-management problem: the monthly reporting unit does not coincide with the weekly constraint unit. Full-horizon modelling includes month boundaries within the model and counts each week as the same seven dates.

The practical requirement is to evaluate conditions continuously across reporting boundaries. Besides the full-horizon approach used here, decomposition carrying history or boundary states can provide that continuity (Kingston, 2021).

## 5.4 Reducing dependence on individual experience

The largest reduction in the task breakdown concerned assigning work and days off, previously the most time-consuming step (Table 6). Computational handling of requirements reduced the need for the responsible staff member to adjust every assignment manually.

Explicit modelling of work rules also provides a basis for turning individual experience into shared specifications. Using draft rosters and exception information as common review materials may support continuity when responsibilities change.

## 5.5 Identifying items for human review

The checking outputs identify items for review in addition to presenting assignments. The weekly checking sheet compares required and actual days off, and the meeting-transfer sheet allows staff to trace date changes applied before solving.

Providing both sheets with the roster allows staff to inspect shortages and applied settings and revise inputs in line with operational schedules. This links the model's scope and exception handling to information available in the workflow.

## 5.6 Threats to validity and scope

### 5.6.1 Construct validity

Preparation time concerns the entire monthly task, whereas solving time concerns a 245-day computation. Day-off shortages and excesses are indicators within defined counting scopes; they do not directly measure operational quality or compliance with every other constraint. Standard remarks on the checking sheet are not a rigorous analysis of infeasibility causes. Learning time and quality after staff reassignment, and the separate effect of the checking sheets, were not measured. Quantitative evaluation of reduced dependence on individual experience remains future work.

### 5.6.2 Internal validity

Working times are operational estimates; they were not obtained through controlled timing with matched months and experience levels. The reduction cannot therefore be attributed separately to individual features. The meeting-day change is a computational scenario, not a report of operational use of a roster with day-off shortages. Neither uniqueness of the affected staff combination nor causal constraints was investigated. Monthly partitioning also changes consecutive-day constraints and full-horizon evaluations, so the observed difference cannot be attributed to a particular objective term.

### 5.6.3 External validity

The evaluation covers 14 staff in the roster studied at one Water Purification Office. The 245 days are a planning horizon, not a post-deployment observation period. Additional investigation is needed for part-time staff, balancing among staff with substantially different active periods, conditions outside the horizon and application to other departments. Results for this input and computing environment should not be directly generalised to other staffing levels, computers or operational effects.

### 5.6.4 Reliability and reproducibility

Constraints, objective terms, weights and stopping conditions are specified, but reproduction also requires corresponding CSV inputs, runtime options, library versions, computing environment and solver logs. We did not study timing distributions over repeated runs or compare other solvers, solving the relaxed model from the outset, or decomposition with history. Changes to work rules, staffing levels, planning horizons, objective terms or weights require model maintenance and renewed checks of weights and solutions. The system can solve any horizon specified by start and end dates. Operational decisions about solving frequency, the horizon of each run and the handling of confirmed assignments are left to the staff and are outside the scope of this evaluation.

# 6 CONCLUSIONS

## 6.1 Main findings

We developed a rostering system combining an integer optimisation model that distinguishes roster states from effective work, full-horizon solving, reference-roster generation when needed, and checking outputs. The system entered routine use at the Water Purification Office of Suita City Waterworks Bureau in May 2026. Operational estimates indicate a reduction in overall work per monthly roster from 12–18 h to approximately 2–3 h.

In experiments evaluating the supporting method, the baseline for 14 staff and 245 days was solved to optimality in approximately 177 s. After the normal model proved infeasible under a changed meeting day, Stage 2 produced a reference roster with five staff-days of shortage. Independent monthly solving without history produced 16 staff–week cases of excess days off across two month-boundary weeks.

The case demonstrates how operational rule modelling, adjustment support and checking information can be integrated into a single workflow, together with the operational effects indicated by the working-time estimates.

## 6.2 Lessons for practice

The application yielded the following lessons:

- Map weekly days off, fixed holiday days off, non-working events and effective work to explicit variables and counting rules.

- Evaluate weekly and consecutive-day conditions across month boundaries through full-horizon modelling or carried assignment history.
- Provide reference-roster generation according to solver status when the normal model yields no feasible solution.
- Present shortages and applied settings with the roster to support staff judgement.
- Distinguish computational experiments, changed-input scenarios and operational estimates in the evaluation.

### 6.3 Future work

Future work will examine other departments and inputs including part-time staff, continuity after staff reassignment, model maintenance and improved presentation of checking information. Larger instances will require reassessment of weights and comparison with approaches such as decomposition with carried history.

## FUNDING

System development was funded internally by TransRecog Co., Ltd. The research and manuscript preparation were partially supported by JSPS KAKENHI Grant Number 25K15007. No funding was received from Suita City.

## DATA AND CODE AVAILABILITY

The datasets generated and/or analysed during the current study are not publicly available due to privacy restrictions but are available from the corresponding author on reasonable request.

## ACKNOWLEDGEMENTS

The authors thank the staff of the Water Purification Office of Suita City Waterworks Bureau for their cooperation in defining operational requirements, providing operational feedback, testing the system and supporting its deployment in routine operations.